\documentclass{crm-article}

\usepackage{amsfonts}
\usepackage{amssymb}
\usepackage{amsmath}
\usepackage{mathpazo}
\usepackage[mathpazo]{flexisym}
\usepackage{comment}
\usepackage{todonotes}
\usepackage[english]{babel}
\usepackage[top=2cm,bottom=2cm,right=2cm,left=2cm]{geometry}
\usepackage{secdot}
\usepackage{hyperref}
\usepackage{url}
\usepackage{mathtools}
\usepackage{algorithm}
\usepackage{algorithmic}
\usepackage{cite}
\usepackage{comment}

\let\argmin\relax
\let\argmax\relax
\let\spn\relax

\DeclareMathOperator{\spn}{span}
\DeclareMathOperator*{\argmax}{arg\,max}
\DeclareMathOperator*{\argmin}{arg\,min}

\newcommand{\la}{\langle}
\newcommand{\ra}{\rangle}

\def\R{\mathbb{R}}

\renewcommand{\leq}{\leqslant}
\renewcommand{\geq}{\geqslant}

\begin{document}

%\englishpaper % раскомментировать в том случае, если текст статьи на английском языке

%\norussian % раскомментировать, если нет метаданных на русском

%\affiliationnoref % раскомментировать, если автор один или все авторы из одной организации

%\emailnoref % раскомментировать в том случае, если автор единственный

%\year=2018 % current year by default
\journalVol{10}
\journalNo{1} %выпуска
\setcounter{page}{1}

% раздел журнала
\journalSection{Математические основы и численные методы моделирования}
\journalSectionEn{Mathematical modeling and numerical simulation}

% дата получения
\journalReceived{01.06.2016.}
%\journalReviewed{01.06.2016.}
%принято к публикации
\journalAccepted{01.06.2016.}

\UDC{519.8}
\title{Ускоренная альтернативная минимизация и адаптивность к сильной выпуклости}
\titleeng{Accelerated alternating minimization and adaptability to strong convexity}
\thanks{Исследование выполнено за счет гранта Российского научного фонда (проект № 18-71-10108)}
\thankseng{This research was funded by Russian Science Foundation (project 18-71-10108).}

%автор - в формате \author{\firstname{И.\,И.}~\surname{Иванов}}
\author{\firstname{Н.\,К.}~\surname{Тупица}}
%автор - в формате \authorfull{Имя Отчество Фамилия}
\authorfull{Назарий Константинович Тупица}
% автор на англ. - в формате \authoreng{\firstname{I.\,I.}~\surname{Ivanov}}
\authoreng{\firstname{N.\,K.}~\surname{Tupitsa}}
%автор на англ. - в формате \authorfull{Firstname M. Surname}
\authorfulleng{Nazarii K. Tupitsa}
%вписать свою электронную почту
%\email{author@example.com}
%организация - в формате \affiliation{Московский государственный университет,\protect\\ Россия, 141700, г. Москва, ул. Университетская, д. 9}
\affiliation{Московский физико-технический институт,\protect\\ Россия, Долгопрудный}
%организация - в формате \affiliationeng{Moscow State Institute University, 9 University street, Moscow, 141700, Russia}
\affiliationeng{Moscow Institute of Physics and Technology,\protect\\ Dolgoprudny, Russia}

% повторите блок для каждого автора;
% если авторов несколько, и автоматическая расстановка сносок от фамилий 
% к организациям приводит к неправильным результатам, укажите правильный 
% вариант в квадраных скобках

\affiliation[1]{Институт проблем передачи и обработки информации,\protect\\ Россия, Москва}
\affiliationeng{Institute for Information Transmission Problems RAS,\protect\\ Moscow, Russia}

\affiliation[1]{Национальный исследовательский университет «Высшая школа экономики»,\protect\\ Россия, Москва}
\affiliationeng{National Research University Higher School of Economics,\protect\\ Moscow, Russia}

\begin{abstract}
В первой части работы приводится ускоренный метод (AGM) первого порядка для выпуклых функций с $L$-липшицевый градиентом, способный автоматически распознавать задачу для которой выполнено условие Поляка-Лоясиевича с константой $\mu$ или являющуюся сильно выпуклой с константой $\mu$. В этих случаях метод демонстрирует линейную сходимость с фактором $\frac{\mu}{L}$, в случае если константа $\mu$ неизвестна. Если значение $\mu$ известно, то линейная сходимоть гарантируется с фактором $\sqrt{\frac{\mu}{L}}$  Если условие не выполнено или задача не является сильно выпуклой, то метод сходится со скоростью $O(1/k^2)$.

Во второй части представлена модификация метода AGM для решения задач, допускающих альтернатвную минимизацию (Alternating AGM). Доказываются аналогичные режимы сходимости.

Таким образом, представлены ускоренные методы демонстрирущие режим линеной сходимости, а именно имеют скорость сходимость $O\left(\min\left\lbrace\frac{LR^2}{k^2},\left(1-{\frac{\mu}{L}}\right)^{(k-1)}LR^2\right\rbrace\right)$, если выполнено условие Поляка-Лоясиевича или являющуюся сильно выпуклой при неизвестной константе $\mu$. Если же условие Поляка-Лоясиевича не выполнено, то сходимость $O(1/k^2)$, но алгоритм при этом не требует никаких изменений.

Также рассматривается адаптивная каталист-оболчка, которая позволяет увеличивать скорость сходимости методов до $O(1/k^2)$, и проводится экспериментальное сравнение градиентного метода с адаптивным выбором шага, ускоренного с помощью адаптивной каталист-оболочки с методами AGM, Alternating AGM, APDAGD и алгоритмом Синхорна на примере задачи двойствненной к задаче оптимального транспорта.

Результатом работы является попытка объяснения причины более быстрой работы метода Alternating AGM  по сравнению с адаптивной каталист оболочкой и методом AGM, несмотря на асимптотически одинаковые гарантии скорости сходимости $O(1/k^2)$. А именно возможной адаптивности к сильной выпуклости метода Alternating AGM. Гипотеза  была проверена на квадратичных задачах, на которых Alternating AGM также оказался быстрее.
\end{abstract}


\begin{abstracteng}
In the first part of the paper we consider accelerated first order optimization method for convex functions with $L$-Lipschitz-continuous gradient, that is able to automatically adapts to problems which satisfies Polyak-Łojasiewicz condition or which is strongly convex with the value of parameter equal to $\mu$. In these cases method possesses linear convergence with factor $\frac{\mu}{L}$, if $\mu$ is unknown. If $\mu$ is known, the the method poses linear convergence with the factor $\sqrt{\frac{\mu}{L}}$. If that are not the cases, the method converges with a rate $O(1/k^2)$. 

The second part contains generalization of the method to the problems, that allows alternating minimization and proofs of the same asymptotic convergence rates.

As a result, it is presented accelerated methods, that are poses linear convergence rate, specifically $O\left(\min\left\lbrace\frac{LR^2}{k^2},\left(1-{\frac{\mu}{L}}\right)^{(k-1)}LR^2\right\rbrace\right)$, if the Polyak-Łojasiewicz condition holds or if the problem is strongly convex with the value of parameter equal to $\mu$. If the Polyak-Łojasiewicz condition does not hold or the problem is not strongly convex, methods poses $O(1/k^2)$ convergence rate, but the methods do not need any adjustment.

Also it is considered the approach called Adaptive Catalyst, which allows to increase convergence rate up to $O(1/k^2)$ and also it is provided an experimental comparison of the approach with AGM, Alternating AGM, APDAGD and Sinkhorn's algorithm for the dual problem to the discrete entropically regularized optimal transport problem.

The result of the work is the attempt to explain the reason why Alternating AGM converge faster than AGM or Adaptive Catalyst despite of the asymptotic theoretical rate $O(1/k^2)$. The hypothesis relies on the fact that Alternating AGM adapts to strong convexity. The hypothesis was tested on quadratic functions, on that Alternating AGM also showed faster convergence.
\end{abstracteng}
\keywordeng{convex optimization}
\keywordeng{alternating minimization}

\maketitle

%Раздел обозначается \paragraph, подраздел - \subparagraph (не \section и \subsection)
\paragraph{Введение}
В данной работе рассматривается задача безусловной оптимизации
\begin{equation}
    \label{eq:pr_st}
    \min_{x \in \R^m} f(x), 
\end{equation}
где $f(x)$ выпусклая с $L$-липшицевым градиентом.
При этом, основным предположением является возможность разделить пространство $\R^m$ на $n$ непересекающихся подпространств $ L_i \in \R^m$, s.t. $\cup L_i = \R^m$ и возможность явно минимизировать функцию $f$ на каждом из подпространств, при прочих фиксированных аргументах. 
Более формально, предполагается, что $f$ имеет блочную структуру $f(x) = f(x_1, \dots, x_n)$, и известно явное решение каждой из задач:
\\
$x^*_i = \argmin_{z\in \R^{n_i}} f(x_1, \dots, x_{i-1}, z, x_{i+1}, \dots, x_n)$, $\sum n_i = m$.

В этих предположениях классической и естественной является идея альтернативной минимизации \cite{ortega2000iterative,bertsekas1989parallel}, где функции минимизируется на каждом подпространстве по очереди. Для гладких сильно выпуклых задач при некоторых дополнительных предположениях была получена линейная сходимость в работе \cite{luo1993error}. В работе \cite{beck2015convergence}  изучен алгоритм альтернативной минимизации для двухблочной функции в достаточно общих предположениях. Предположения заключаются в существовании хотя одной гладкой блочной компоненты. Также показано, что негладкость допустима композитного слагаемого, что не влият свойства сходимости. Поскольку среди этих предположений нет сильной выпуклости, получена лишь сублинейная скорость сходимости $O(1/k)$, где $k$ номер итерации. Похожий результат получен для произвольного числа блоков \cite{hong2017iteration,sun2015improved}. В предположениях гладкости и сильной выпуклости в работе \cite{nutini2015coordinate} получена линейная сходимость для произвольного числа блоков. Скорость сходимости пропорциональна $\kappa$ -- эффективному числу обусловленности задачи. В работе \cite{chambolle2017accelerated} представлен ускоренный алгоритм альтернативнйо минимизации для задачи специального вида и для двух блоков, а именно с задачей вида суммы квадратичных функций с композитными членами проксимального вида. Получена скорость сходимости порядка $O(1/k^2)$ для выпуклых постановоки  линейная с экспонентой  $\sqrt{\kappa}$ в сильно выпуклом случае. В работе \cite{diakonikolas2018alternating} рассматрен неускоренный алгоритм альтернативной минимизации и получена скорость сходимости $O(1/k)$в выпуклом случае и линейная сходимость с экспонентой  $\kappa$ в сильновыпуклом случае. Также предложен ускоренный метод для задаче в общей выпуклой постановке со сходимостью $O(1/k^2)$ и предположение о возможности обобщения подхода к рассмотрению на сильно выпуклый случай. В работе \cite{2019arXiv190603622G} рассматривается ускоренный метод бщей выпуклой постановке со сходимостью $O(1/k^2)$ и $\sqrt{\kappa}$ в сильно выпуклом случае. для произвольного числа блоков. Также стоит упомянуть обзорную работу \cite{hong2016unified}.

В первой части работы доказываем адаптивность алгоритма 1 из \cite{nesterov2020primaldual} результат, который, по всей видимости, ранее не был опубликован, несмотря на обширные исследования различных модификаций метода \cite{guminov2019accelerated}. Таким образом получаем метод демонстрирующий скорость сходимости $O(1/k^2)$ в гладком случае и линейную с фактором  ${\kappa}$ в сильно выпуклом случае, причем при неизвестном значении параметра сильной выпуклости. Также обобщаем этот результат на случай альтернативной минимизации. Более того линейную сходимость получить при более слабом условии Поляка-Лоясиевича \cite{polyak1987introduction}.

Далее приводим обобщение ускоренного алгоритма альтернативной минимизации из работы \cite{2019arXiv190603622G} (алгоритм 1) на сильно выпуклый случай.

Также рассматривается другой подход к ускорению альтернативной минимизации, основанный на работе \cite{2019arXiv191111271I}, и проводится его сравнение с алгоритмом 1 из \cite{nesterov2020primaldual} и неускоренным методом альтернативной минимизации на примере задачи двойственной к задаче оптимального транспорта.

\paragraph{Адаптивность к сильной выпуклости ускоренного градиентного метода}

Рассмотрим алгоритм 3 из работы \cite{nesterov2020primaldual}. Для его работы требует указать значение параметра сильной выпуклости оптимизируемой функции, тогда его сложность описывается выражением
\[
    f(x_{k})-f(x_*)\leqslant\min\left\lbrace\frac{2LR^2}{k^2},\left(1-\sqrt{\frac{\mu}{L}}\right)^{(k-1)}LR^2\right\rbrace.
\]
Если же значение параметра сильной выпоклости неизвестно, то можно запустить алгоритм с $\mu=0$. В таком случае алгоритм 3 из \cite{nesterov2020primaldual} будет в точностью представлять собой алгоритм 1 из той же работы, который приведен ниже. В этом случае авторы работы \cite{nesterov2020primaldual} гаранируют лишь сходимость со скорость
\[
    f(x_{k})-f(x_*) \leqslant \frac{2LR^2}{k^2}
\]

Ниже приведен алгоритм 1 из \cite{nesterov2020primaldual}
\floatname{algorithm}{Алгоритм}

\begin{algorithm}[!h]
\caption{Accelerated Gradient Method with Small-Dimensional Relaxation (AGMsDR)}
\label{AGMsDR}
\begin{algorithmic}[1]
%\REQUIRE $x^0 = v^0$, $L$
\ENSURE $x^k$
\STATE Полагаем $k = 0$, $A_0=0$, $x^0 = v^0$, $\psi_0(x) = V(x,x^0)$
\FOR{$k \geq 0$}
\STATE \begin{equation}
    \label{eq:beta_k_y_k_def}
    \beta_k = \arg\min_{\beta \in \left[0, 1 \right]} f\left(v^k + \beta (x^k - v^k)\right), \quad y^k = v^k + \beta_k (x^k - v^k).
\end{equation}
\STATE 
Вариант а), значение $L$ известно, 
\begin{equation}
\label{eq:xkp1_opt_a}
    x^{k+1} = \arg\min_{x \in E} \left\{ f(y^k) + \langle \nabla f(y^k), x - y^k \rangle + \frac{L}{2} \|x - y^k \|_2^2 \right\}.
\end{equation}
Найти $a_{k+1}$ из уравнения $\frac{a_{k+1}^2}{A_{k} + a_{k+1}} = \frac{1}{L}$.  \\
Вариант б), 
\begin{equation}
\label{eq:xkp1_opt_b}
h_{k+1} = \arg\min_{h \ge 0} f\left(y^k - h(\nabla f(y^k))^{\#}\right), \quad x^{k+1} = y^  k- h_{k+1}(\nabla f(y^k))^{\#}.
\end{equation}
Найти $a_{k+1}$ из уравнения $f(y^k) - \frac{a_{k+1}^2}{2(A_{k} + a_{k+1})} \|\nabla f(y^k) \|_2^2 = f(x^{k+1})$.
\STATE  Полагаем $A_{k+1} = A_{k} + a_{k+1}$.  
\STATE  Полагаем $\psi_{k+1}(x) =  \psi_{k}(x) + a_{k+1}\{f(y^k) + \langle \nabla f(y^k), x - y^k \rangle\}$.
\STATE $v^{k+1} = \arg\min_{x \in E} \psi_{k+1}(x)$
\STATE $k = k + 1$
\ENDFOR
\end{algorithmic}
\end{algorithm}

Следующая лемма объясняет поведение Алгоритма \ref{AGMsDR}, если оптимизируемая функция является сильно выпуклой, но значения параметра неизвестно.
Но перед этим обртатимся к \cite{polyak1987introduction}, где показано, что сильная выпуклость влечет за собой более слабое условие Поляка-Лоясиевича, которое выполняется для более широкого класса задач.
\begin{lem}
    Алгоритм \ref{AGMsDR} автоматически распознает сильно выпуклую задачу или задачу для которой выполнено условие Поляка-Лоясиевича и демонстрирует линейную скорость сходимости:
    \[ f(x^{k+1}) - f(x^*) \leq \Pi_{i=0}^{k-1}  \big(1 - \frac{\mu}{\hat L_i}\big) \cdot (f(x^{0})- f(x^{*})),\]
    где $\hat L_i = \frac{A_{i} + a_{i+1}}{a_{i+1}^2}$ является оценкой локальной константы Липшица градиента функции $L$ на $i$-й итерации.
    \label{non-acc-conv}
\end{lem}

\proof*
    Доказательство леммы основано на доказательстве теоремы 1 из \cite{nesterov2020primaldual}, откуда было взято неравенство
    \begin{multline*}
        f(x^{k+1}) \leq f(y^k) + \langle \nabla f(y^k), x^{k+1} - y^k \rangle + \frac{L}{2} \|x^{k+1} - y^k \|_2^2
        \\
        = \min_{x \in E} \left(f(y^k) + \langle \nabla f(y^k), x - y^k \rangle + \frac{L}{2} \|x - y^k \|_2^2 \right)
        = f(y^k) - \frac{1}{2L} \|\nabla f(y^k)\|_2^2.
    \end{multline*}
    Это неравенство задает условие достаточного убывания
    \begin{equation*}
        f(y^k) - \frac{a_{k+1}^2}{2(A_{k} + a_{k+1})} \| \nabla f(y^k) \|_2^2  = f(x^{k+1}) \geq f(y^{k+1})
    \end{equation*}
    где последнее неравенство получено в силу \eqref{eq:beta_k_y_k_def}, откуда следует, что $f(y^k)\leq f(x^k)$.
    Скомпоновав данный результат с условием Поляка-Лоясиевича
    \begin{equation*}
         \|\nabla f(y^k)\|_2^2  \geq 2\mu\left(F(y^{k})- F(x^*)\right)
    \end{equation*}
    которое является следствием сильной выпуклости, при этом константа в ПЛ-условии совпадает с константой сильной выпуклости,
    получим 
    \[
        \left(f(y^{k+1})- f(x^*)\right) \leq (1 - \frac{\mu a_{k+1}^2}{A_{k} + a_{k+1}}) \left(f(y^{k})- f(x^{*})\right) 
        \leq \Pi_{i=0}^{k}  (1 - \frac{\mu a_{i+1}^2}{A_{i} + a_{i+1}}) (f(x^{0})- f(x^{*}))
    \]
    И оканчательно, воспользуемся условием \eqref{eq:xkp1_opt_b}, которое обеспечивает $f(x^{k+1})\leq f(y^k)$, получим линейную сходимость
    \[ f(x^{k+1}) - f(x^*) \leq \Pi_{i=0}^{k-1}  (1 - \frac{\mu a_{i+1}^2}{A_{i} + a_{i+1}}) (f(x^{0})- f(x^{*}))\]
    Также заметим, что $\frac{A_{i} + a_{i+1}}{a_{i+1}^2}$ являтеся оценкой сверху локальной константы Липшица градиента функции на $i$-ой итерации .
\qed

Заметим также, что для глобальной консанты Липшица градиента функции может выполняться $L \geq \frac{A_{i} + a_{i+1}}{a_{i+1}^2}$,  что означает, что скорость сходимости может оказаться лучше чем с фактором $\kappa=\frac{\mu}{L}$, где $L$ - глобальная константа Липшица градиента. 

Таким образом данный алгоритм демонстрирует нижние оценки сложности на функциях с Липшицевых градиентом и адаптивность к сильно выпуклым функциям, но в последнем случае нижние границы сложности не достигаются.

\paragraph{Альтернативная минимизация}
В этой секции рассмотрим вариацию Алгоритма \ref{AGMsDR} для задач, допускающих альтернативную минимизацию. В работе \cite{2019arXiv190603622G} рассмотрена такая вариация для функций с липшицевым градиентом. Далее приводится модификация алгоритма для решения сильно выпуклых задач, допускающих альтернативную минимиизацию.

Для удобства введем обозначения. Множество $\{1,\ldots, N\}$ векторов $\{e_i\}_{i=1}^N$ ортонормированного базиса разделено на $n$ непересекающихся блоков  $I_k$, $k\in\{1,\ldots,n\}$(см. Введение). Пусть $S_k(x)=x+\spn\{e_i:\ i\in I_k\}$, подпространство, содержащее $x$ построенное на базисных векторах $k$-го блока.

\begin{algorithm}[H]
\caption{Accelerated Alternating Minimization}
\label{AAM-2}

\begin{algorithmic}[1]
   \REQUIRE Starting point $x_0$
    \ENSURE $x^k$
   \STATE Полагаем $A_0=0$, $x^0 = v^0$, $\tau_0 = 1$
   \FOR{$k \geqslant 0$}
	\STATE Полагаем 
	\begin{equation}
	    \beta_k = \argmin_{\beta \in [0,1]} f\Big(x^k + \beta (v^k - x^k)\Big)
	    \label{line-search}
	\end{equation}
	\STATE Полагаем $y^k = x^k + \beta_k (v^k - x^k) $
    \STATE Выбираем $i_k=\argmax_{i\in\{1,\ldots,n\}} \|\nabla_i f(y^k)\|_2^2$
\STATE Полагаем $x^{k+1}=\argmin_{x \in S_{i_k}(y^k)} f(x)$
\STATE Если $L$ известно находим $a_{k+1}$ из уравнения
$\frac{a_{k+1}^2}{(A_{k}+a_{k+1})(\tau_{k}+\mu a_{k+1})} = \frac{1}{Ln}$
\\
Если $L$ не известно находим $a_{k+1}$ из уравнения
\begin{multline}
    f(y^k) - \frac{a_{k+1}^2}{2(A_{k} + a_{k+1})(\tau_{k}+\mu a_{k+1})}\| \nabla f(y^k) \|_2^2 +
    \\
    \frac{\mu \tau_k a_{k+1}}{2(A_{k} + a_{k+1})(\tau_{k}+\mu a_{k+1})}\| v^k - y^k \|_2^2 = f(x^{k+1})
    \label{aam-s.d.}
\end{multline}

\STATE  Полагаем $A_{k+1} = A_{k} + a_{k+1}$, $\tau_{k+1} = \tau_k + \mu a_{k+1}$
%\STATE  Полагаем $\psi_{k+1}(x) =  \psi_{k}(x) + a_{k+1}\{f(y^k) + \langle \nabla f(y^k), x - y^k \rangle\}$
\STATE Полагаем $v^{k+1}=\argmin_{x\in\mathbb{R}^N} \psi_{k+1}(x)$ 
%\STATE $k = k + 1$
\ENDFOR
\end{algorithmic}

\end{algorithm}

Введем вспомогательную последовательность функций \[\psi_0(x)=\frac{1}{2}\|x-x^0\|_2^2,\]

\[\psi_{k+1}(x) = \psi_{k}(x) + a_{k+1}\{f(y^k) + \langle \nabla f(y^k), x - y^k \rangle\ + \frac{\mu}{2}\| x - y^k \|_2^2\}.\]

Используя следующее обозначении
\[ l_k(x) = \sum_{i=0}^k a_{i+1} \{ f(y^i) + \langle \nabla f(y^i), x - y^i \rangle\ + \frac{\mu}{2}\| x - y^i \|_2^2 \} \]
получим нереккурентное представление
\[
\psi_{k+1}(x) = \psi_{0}(x) + l_k(x)
\]
Заметим, что $\psi_k(x)$ is $\tau_k$ является сильно выпуклой функцией со значением параметра $\tau_k = 1 + \mu\sum_{i=0}^k a_i = 1 + \mu A_k$.
\begin{lem} 
\label{AAM-2_Ak_rate}
После $k$ итераций Алгоритма \ref{AAM-2} выполняется
\begin{equation}
    \label{eq:main_recurrence}
    A_{k}f(x^{k}) \leqslant \min_{x \in \mathbb{R}^N} \psi_{k}(x) = \psi_{k}(v^{k}).
\end{equation}
Более того, если функция с $L$ липшицевым градиентом и $\mu\geq0$ сильно выпукла, то
\[A_k \geqslant \max \left\{\frac{k^2}{4Ln}, \frac{1}{nL}\left(1 - \sqrt{\frac{\mu}{nL}}\right)^{-k-1}\right\},\]
где $n$ количество блоков, по которым допускается явная минимизация.
\end{lem}

\proof*
Докажем \eqref{eq:main_recurrence} индукцией по $k$. При $k=0$ неравенство верно. Предположим
\[A_{k}f(x^{k}) \leqslant \min_{x \in \mathbb{R}^N} \psi_{k}(x) = \psi_{k}(v^{k}).\]
Тогда
\begin{multline*}
    \psi_{k+1}(v^{k+1}) = \min_{x \in \mathbb{R}^N} \Bigg\{ \psi_{k}(x) + a_{k+1}\{f(y^k) + \langle \nabla f(y^k), x - y^k \rangle + \frac{\mu}{2}\| x - y^k \|_2^2\} \Bigg\}
    \\
    \geqslant \min_{x \in \mathbb{R}^N} \Bigg\{ \psi_{k}(v^k) + \frac{\tau_k}{2}\| x - v^k \|_2^2 + a_{k+1}\{f(y^k) + \langle \nabla f(y^k), x - y^k \rangle 
    \\
    \shoveright{+ \frac{\mu}{2}\| x - y^k \|_2^2\} \Bigg\}}
    \geqslant \min_{x \in \mathbb{R}^N} \Bigg\{ A_{k}f(x^k) + \frac{\tau_k}{2}\| x - v^k \|_2^2 + a_{k+1}\{f(y^k) + \langle \nabla f(y^k), x - y^k \rangle 
    + \frac{\mu}{2}\| x - y^k \|_2^2\} \Bigg\},
\end{multline*}
где была исполльзовано, что $\psi_{k}$ сильно выпуклая и имеет минимум в $v^k$ и то, что $f(y^k)\leqslant f(x^k)$. 

Условия оптимальности для $\min\limits_{\beta\in [0,1]} f\left(x^k + \beta (v^k - x^k)\right)$ гарантируют одно из
\begin{enumerate}
    \item $\beta_k = 1$, $\langle \nabla f(y^k),x^k - v^k \rangle \geqslant 0$, $y^k = v^k$;
    \item $\beta_k \in (0,1)$ and $\langle \nabla f(y^k),x^k - v^k \rangle = 0$, $y^k = v^k + \beta_k (x^k - v^k)$;
    \item $\beta_k = 0$ and $\langle \nabla f(y^k),x^k - v^k \rangle \leqslant 0$, $y^k = x^k$ .
\end{enumerate}
Во всех случаях будет верно  $\langle \nabla f(y^k), v^k - y^k \rangle \geqslant 0$.

Таким образом
\begin{multline*}
    \psi_{k+1}(v^{k+1}) \geqslant  \min_{x \in \mathbb{R}^N} \Big\{ A_{k}f(y^k) + \frac{\tau_k}{2}\| x - v^k \|_2^2 + a_{k+1}\{f(y^k) + \langle \nabla f(y^k), x - y^k \rangle
    + \frac{\mu}{2}\| x - y^k \|_2^2\} \Big\}
\end{multline*}
Решение задачи выше находится явно:
\[
x=\frac{1}{\tau_{k+1}}(\tau_k v^k + \mu a_{k+1}y^k - a_{k+1} \nabla f(y^k)),
\]
подставим его и используем то, что $\langle \nabla f(y^k), v^k - y^k \rangle \geqslant 0$, получим
\begin{align*}
    \psi_{k+1}(v^{k+1}) &\geqslant A_{k+1}f(y^k) - \frac{a_{k+1}^2}{2\tau_{k+1}}\| \nabla f(y^k) \|_2^2 + \frac{\mu \tau_k a_{k+1}}{2\tau_{k+1}}\| v^k - y^k \|_2^2.
\end{align*}

Далее покажем, что
\begin{align*}
    A_{k+1}f(y^k) - \frac{a_{k+1}^2}{2\tau_{k+1}}\| \nabla f(y^k) \|_2^2 + \frac{\mu \tau_k a_{k+1}}{2\tau_{k+1}}\| v^k - y^k \|_2^2 &\geqslant A_{k+1}f(x^{k+1})
\end{align*}
что завершит индукционный переход

Для этого, принимая во внимание, что $f$ имеет $L$-липшицев градиент, получим
$\forall i$
\[
f(y^k)-\frac{1}{2L}\|\nabla_i f(y^k)\|_2^2\geqslant f(x_i^{k+1}),
\]
где $x_i^{k+1}=\argmin_{x\in S_{i}} f(x)$. Так как $i_k=\argmax_{i} \|\nabla_i f(y^k)\|_2^2$, $$\|\nabla_{i_k} f(y^k)\|_2^2\geqslant \frac{1}{n}\|\nabla f(y^k)\|_2^2$$ 
и
$$f(y^k)-\frac{1}{2Ln}\|\nabla f(y^k)\|_2^2\geqslant f(y^k)-\frac{1}{2L}\|\nabla_{i_k} f(y^k)\|_2^2\geqslant f(x^{k+1}),$$
Выбор $a_{k+1}$ такого, что $\frac{a_{k+1}^2}{2A_{k+1}\tau_{k+1}}\geqslant \frac{1}{2Ln}$ означает 
%\frac{a_{k+1}^2}{2A_{k+1}\tau_{k+1}}\geqslant \frac{1}{2Ln}
\begin{multline*}
    A_{k+1}f(y^k) - \frac{a_{k+1}^2}{2\tau_{k+1}}\| \nabla f(y^k) \|_2^2 + \frac{\mu \tau_k a_{k+1}}{2\tau_{k+1}}\| v^k - y^k \|_2^2 
    \\
    \geqslant A_{k+1}f(y^k) - \frac{a_{k+1}^2}{2\tau_{k+1}}\| \nabla f(y^k) \|_2^2  
    \geqslant A_{k+1}f(y^k) - \frac{A_{k+1}}{2Ln} \| \nabla f(y^k) \|_2^2     
    \geqslant A_{k+1}f(x^{k+1}),
\end{multline*}
что завершает доказательство шага индукции.

Преобразовав выражение для поиска $a_{k+1}$ получим $\frac{a_{k+1}^2}{(A_{k}+a_{k+1})(\tau_{k}+\mu a_{k+1})}\geqslant \frac{1}{Ln}$.

Теперь оценим скорость роста последовательности $A_k$. $\tau_k = 1 + \mu\sum_{i=0}^k a_i = 1 + \mu A_k$. $\frac{a_{k+1}^2}{2A_{k+1}\tau_{k+1}}\geqslant \frac{1}{2Ln}$
\begin{align*}
	a^2_{k} & \geqslant 
	\frac{A_k \tau_k}{{nL}} = \frac{{A_k + \mu A_k^2}}{nL}
\end{align*}

\begin{align}
	a_{k} & \geqslant 
	\frac{1}{\sqrt{nL}} \sqrt{A_k + \mu A_k^2}
	\geqslant 
	\sqrt{\frac{\mu}{2Ln}}A_k
\end{align}

\begin{equation*}
    \sqrt{A_i} - \sqrt{A_{i-1}}
    \geqslant \frac{A_i-A_{i-1}}{\sqrt{A_i} + \sqrt{A_{i-1}}} 
    \geqslant \frac{a_i}{2\sqrt{A_i}}
    \geqslant \frac{\sqrt{1+\mu A_i}}{2\sqrt{Ln}}
\end{equation*}

Суммируя по $i=1,\dots,k$ получим \[A_k \geqslant \frac{k^2}{4Ln}\]

Также \[A_{k+1} = A_{k} + a_{k+1} \geqslant A_k + \sqrt{\frac{\mu}{nL}}A_{k+1},\] что означает \[A_{k+1} \geqslant \left(1 - \sqrt{\frac{\mu}{nL}}\right)^{-1}A_{k}\]

Остается оценить $A_1$:
\[A_1 = \frac{a_1^2}{A_1} \geqslant  \frac{a_1^2} {(1 + \mu A_1)A1} \geqslant  \frac{a_1^2}{A_1 \tau_1} \geqslant \frac{1}{nL}\]

реккурсивно применяя полученные оценки, приходим к утверждению леммы:
\[A_k \geqslant \max \left\{\frac{k^2}{4Ln}, \frac{1}{nL}\left(1 - \sqrt{\frac{\mu}{nL}}\right)^{-k+1}\right\}\]
\qed

\begin{teo} 
После $k$ итераций Алгоритма \ref{AAM-2} выполняется

\begin{equation}
    \label{eq:main_result}
    f(x^k) - f(x_*) \leqslant nLR^2\min \left\{ \frac{4}{k^2}, \left(1 - \sqrt{\frac{\mu}{nL}}\right)^{k-1}\right\}
\end{equation}
\end{teo}

\proof*
From the convexity of $f(x)$ we have
\[
    l_k(x_*) = \sum_{i=0}^{k} a_{i+1}(f(y^i)+\langle\nabla f(y^i),x_*-y^i\rangle + \frac{\mu}{2}\| x_* - y^i \|_2^2) 
\leqslant A_{k+1} f(x_*).
\]
По Лемме \ref{AAM-2_Ak_rate}  
\begin{multline*}
    A_k f(x^k)\leqslant\psi_{k}(v^k)
    \leqslant\psi_k(x_*)=\frac{1}{2}\|x_*-x^0\|_2^2 \notag \\
    +\sum_{i=0}^{k-1} a_{i+1}(f(y^i)+\langle\nabla f(y^i),x_*-y^i\rangle + \frac{\mu}{2}\| x_* - y^i \|_2^2) 
    \leqslant 
     A_k f(x_*)+\frac{1}{2}\|x_*-x^0\|_2^2
\end{multline*}

И окончательно
\begin{equation*}
    f(x^k) - f(x_*) \leqslant \frac{R^2}{2A_k} \leqslant nLR^2\min \left\{ \frac{4}{k^2}, \left(1 - \sqrt{\frac{\mu}{nL}}\right)^{k-1}\right\}.
\end{equation*}
\qed

\begin{lem}
    Алгоритм \ref{AAM-2} автоматически распознает сильно выпуклую задачу или задачу для которой выполнено условие Поляка-Лоясиевича и демонстрирует линейную скорость сходимости:
    \[ f(x^{k+1}) - f(x^*) \leq \Pi_{i=0}^{k-1}  \big(1 - \frac{\mu}{\hat L_i}\big) \cdot (f(x^{0})- f(x^{*})),\]
    где $\hat L_i = \frac{A_{i} + a_{i+1}}{a_{i+1}^2}$ является оценкой локальной константы Липшица градиента функции $L$ на $i$-й итерации.
\end{lem}
\proof*
Доказательство повторяет доказательство Леммы \ref{non-acc-conv}.
\qed

\paragraph{Экспериментальное сравнение}
Большое количество исследований посвящено решению задач оптимального транспорта \cite{cuturi2013sinkhorn, pmlr-v80-dvurechensky18a, 2019arXiv190603622G}, поиску барицентров Вассерштейна \cite{kroshnin2019complexity, dvinskikh2019primal, uribe2018distributed, dvurechensky2018decentralize, 2020arXiv200204783L}, а также много-маргинального оптимального транспорта \cite{2019arXiv191000152L, 2020arXiv200402294T} двойственными методами. Все эти задачи допускают альтернативную минимизацию. Далее проводится экспериментальное сравнение представленных алгоритмов на примере задачи двойственной к энтропийно регуляризованной задаче дискретного оптимального транспорта (ЭОТ) \cite{2019arXiv190603622G}.

Как известно \cite{2019arXiv190603622G}, задача ЭОТ выглядит следующим образом:
\begin{equation}
    \label{OT_dual}
    f(u,v)=\gamma(\ln\left(\mathbf{1}^TB(u,v)\mathbf{1}\right)-\la u,r \ra - \la v,c \ra)\to \min_{u,v\in\mathbb{R}^N} ,
\end{equation}
где
%$\vp(u,v)=\gamma(\ln\left(\mathbf{1}^TB(u,v)\mathbf{1}\right)-\la u,r \ra - \la v,c \ra)$,
$[B(u,v)]^{ij} = \exp \left(u^{i}+v^{j} -\frac{C^{i j}}{\gamma}\right)$,  $B,C \in \R ^{N \times N}_+$, $\gamma \in \R_+$.
Переменные в этой задаче естественным образом разделяются на два блока, и при фиксированных переменных одного блока удается явно выписать решение условий оптимальности по другому блоку. Таким образом поочередно обновляя переменные и получается алгоритм Синхорна.

В работе \cite{pmlr-v80-dvurechensky18a} представлен алгоритм APDAGD, который показал более быструю сходимость по сравнению с алгоритмом Синхорна.

В работе \cite{guminov2019accelerated} Алгоритм \ref{AAM-2}, примененный к этой задаче, показал себя наиболее стабильным и быстрым в сравнении с другими алгоритмами образом.

\begin{figure}[H]

\centering
\includegraphics{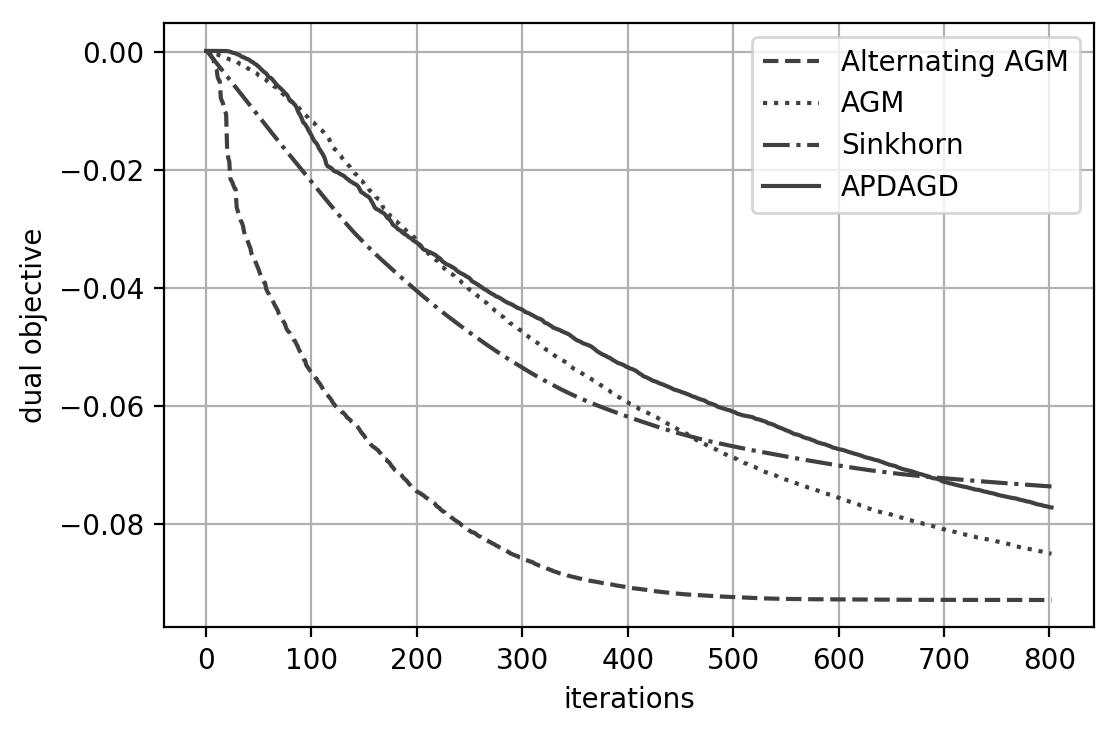}
\smallskip
\caption{Сходимость алгоритма Синхорна и методов AGM, Alternating AGM и APDAGD, примененных к ЭОТ}
\label{agm}
\end{figure}

На рисунке \ref{agm} приводится эспериментальное сранение Алгоритмов \ref{AGMsDR} и \ref{AAM-2}, а также алгоритма Синхорна и алгоритма APDAGD на примере задачи ЭОТ.

Следует заметить, что теоретические оценки сложности для Alternating AGM показыают, что для достижения заданной точности ему требуется в $\sqrt{n}$ раз большее число итераций, чем для AGM, где $n$ - число блоков по которым возможна альтернативная минимизация, причем для ЭОТ $n=2$. Экспериментальное сравнение дает противополжные результаты.
Одной из целей данной работы является попытка объяснения этого явления. Предполагается, что такое поведение связано с адаптивностью методов Alternating AGM и AGM к сильной выпуклости. Как известно, задача ЭОТ не является сильно-выпуклой так как значение целевой функции инвариантно на прямых, параллельных вектору $(\mathbf{1},-\mathbf{1})$. Поэтому схожая сходимость методов AGM и APDAGD является ожидаемой, т.к. в не сильно выпуклом случае скорость сходимости имеет порядок $O(1/k^2)$. Alternating AGM также имеет сходимость порядка $O(1/k^2)$ для не сильно выпуклых функций, но, вероятно, распознает сильную  выпуклость на подпространствах переменных $u$ и $v$, ортогональных вектору $(\mathbf{1},-\mathbf{1})$, что может приводить к  наблюдаемой более быстрой сходимости.

Данная гипотеза была также проверена на квадратичных задачах 
\begin{equation}
    \min_z f(z) = \|Wz-b\|_2^2 .
    \label{qprob}
\end{equation}
$f$ является сильно-выпуклой с константой $\mu = \sqrt{\lambda_{\min}(W^T W)}$. Матрица $W$ является симметричной. 

Последняя задача может быть решена с помощью Алгоритма \ref{AGMsDR}.

Построим эквивалентную задачу, допускающую альтернативную минимизацию путем разделения вектора $z$ на два блока одинакового размера:
\[
    z = 
    \begin{pmatrix}
        x
        \\
        y
    \end{pmatrix}.
\]
Также разделим матрицу $W$ на 4 блока одинакового размера
\[
    W = 
    \begin{pmatrix}
        A B
        \\
        C D
    \end{pmatrix}
\]
и вектор $b$ на два
\[
    b = 
    \begin{pmatrix}
        c
        \\
        d
    \end{pmatrix}.
\]

Тогда задача (\ref{qprob}) будет эквивалентна задаче
\begin{equation}
    \min_{x,y} \|Ax + By - c\|_2^2 + \| Cx + Dy -d\|_2^2.
    \label{qprob2}
\end{equation}
Последняя задача допускает альтернативную минимизацию
\begin{align*}
    x^{k+1} &= (A^T A + C^T C)^{-1}\big[A^T(c-By^k) + C^T(d-Dy^k)\big]
    \\
    y^{k+1} &= (B^T B + D^T D)^{-1}\big[B^T(c-Ax^k) + D^T(d-Cx^k)\big]
\end{align*}
и может быть решена с помощью Алгоритма \ref{AAM-2}.

\begin{figure}[H]
\centering
\includegraphics[width=\textwidth]{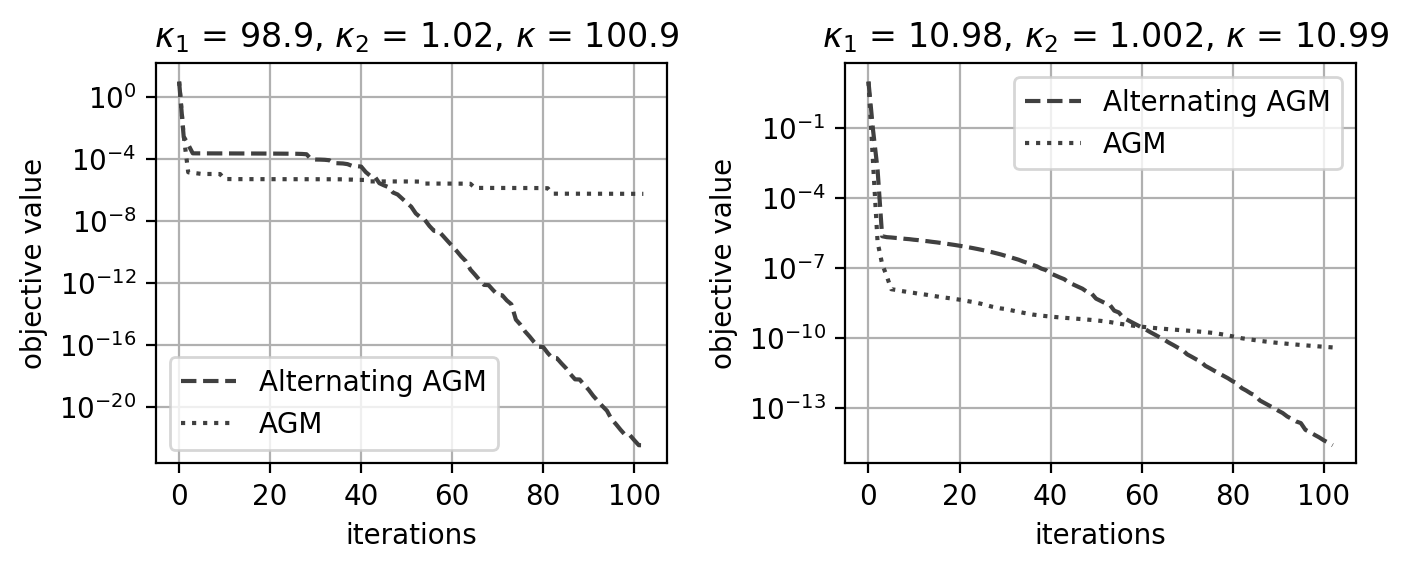}
\smallskip
\caption{Сходимость алгоритмов AGM и Alternating AGM, примененных к задачам \eqref{qprob} и \eqref{qprob2} соответственно}
\label{quadr-1}
\end{figure}

\begin{figure}[H]
\centering
\includegraphics[width=\textwidth]{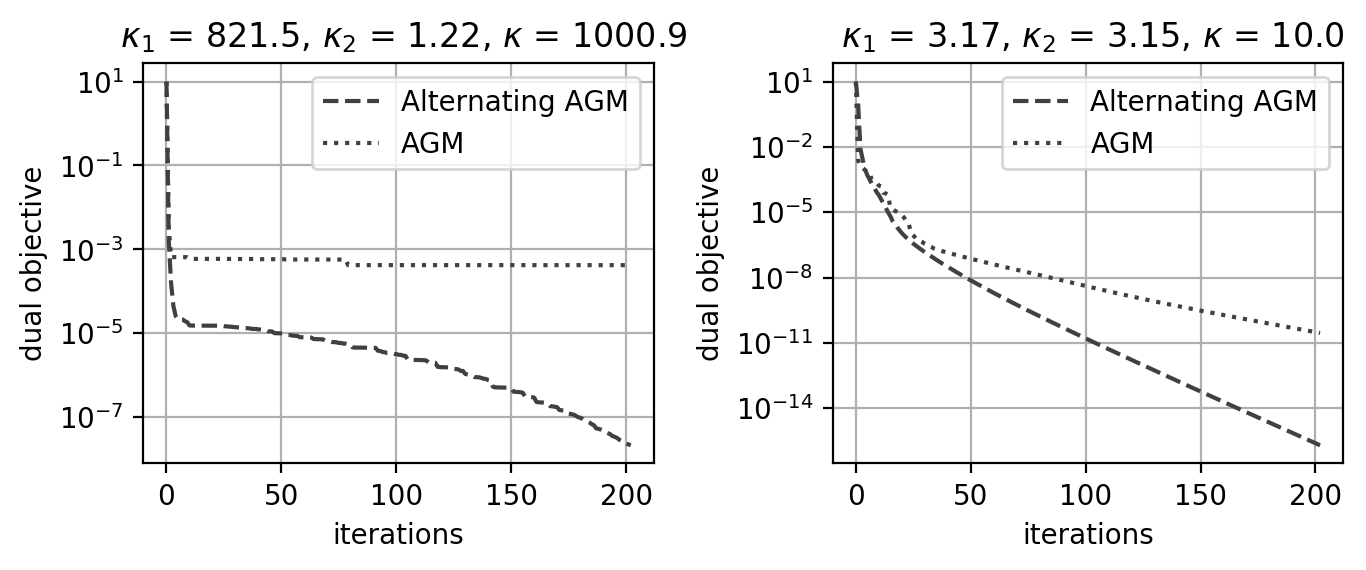}
\smallskip
\caption{Сходимость алгоритмов AGM и Alternating AGM, примененных к задачам \eqref{qprob} и \eqref{qprob2} соответственно}
\label{quadr-2}
\end{figure}

Результаты сравнения алгоритмов AGM и Alternating AGM, запущенных с $\mu=0$ представленны на рисунках \ref{quadr-1} и \ref{quadr-2} для различных чисел обусловленности матрицы. $\kappa, \kappa_1, \kappa_2$ -- числа обусловленности матриц $W$, $A$ и $D$ соответственно. По всей видимости, более быстрая сходимость связана с тем, что один из блоков (или оба блока) обусловленны лучше, чем вся задача.

\paragraph{Адаптивная каталист оболочка}
В работе \cite{2019arXiv191111271I} рассматривается так называемая адаптивная каталист-оболочка, позволяющая ускорять методы.
Для удобства приводим обозначения и алгоритм
\begin{eqnarray*}
    F_{L,x}(y)& = & f(y) + \tfrac{L}{2}||y-x||_2^2,\\
\end{eqnarray*}
$f$ - минимизируемая функция.
\begin{algorithm}[ht!]
\caption{Adaptive Catalyst}
\label{alg1}
\begin{algorithmic}[1]
\REQUIRE Начальная точка $x^0$, оценка $L_0 > 0$, параметры $\alpha > \beta > \gamma > 0$ и метод $\mathcal{M}$.
\STATE Полагаем $y^0=z^0=x^0$
%\Repeat
\FOR{$k = 0, 1, \ldots,N-1$}
    \STATE $L_{k+1} =  \beta \cdot\min\left\{ \alpha L_{k},L_u\right\}$ 
    \STATE $t = 0$
    \REPEAT
    \STATE{$t: = t+1$}
     \STATE{$L_{k+1}:=  \max\left\{L_{k+1}/ \beta,L_d\right\}$}
        \STATE Вычисляем
           \begin{eqnarray*}
            a_{k+1} &=& \tfrac{1/L_{k+1} + \sqrt{1/L^2_{k+1} + 4 A_k/L_{k+1}}}{2},\\
            A_{k+1} &=& A_k + a_{k+1},\\
            x^{k+1}& =& \tfrac{A_{k}}{A_{k+1}} y^{k} + \tfrac{a_{k+1}}{A_{k+1}} z^{k}.
              \end{eqnarray*}
        \STATE Вычисляем приближенное решение следующей задачи спомощью вспомогательного неускоренного метода $\mathcal{M}$
        $$ y^{k+1} \approx  \argmin_{y}  F_{L, x^{k+1}} (y)   $$
         \STATE Запуская метод $\mathcal{M}$ из точки $x^{k+1}$ и ожидая на выходе точку $y^{k+1}$ делаем  $N_t$ итераций и проверяем адаптивный критерий остановки:
        \begin{equation}
            \|\nabla F_{L, x^{k+1}}  (y^{k+1})\|_2 \leq \tfrac{L_{k+1}}{2} \|y^{k+1} - x^{k+1}\|_2.  % (Monteiro-Svaiter condition)
            \label{MS-cond}
        \end{equation} 
    \UNTIL {$t>1$ and $N_{t} \ge \gamma \cdot N_{t-1}$ or $L_{k+1} = L_d$}
    \STATE $z^{k+1}=z^{k}-a_{k+1} \nabla f\left(y^{k+1}\right)$
%\Until {$f\left(y^{k+1}\right) - f\left(x^*\right) \le \varepsilon$}
% \Ensure{$y^{k+1}$}
\ENDFOR
\STATE {\bf Output:} $y^{N}$

\end{algorithmic}
\end{algorithm}

Данный алгоритм был экспериментально изучен на задаче ALS и оказался успешным. Задача и подход к решению описаны в \cite{4781121}.

\begin{figure}[H]
\centering
\includegraphics{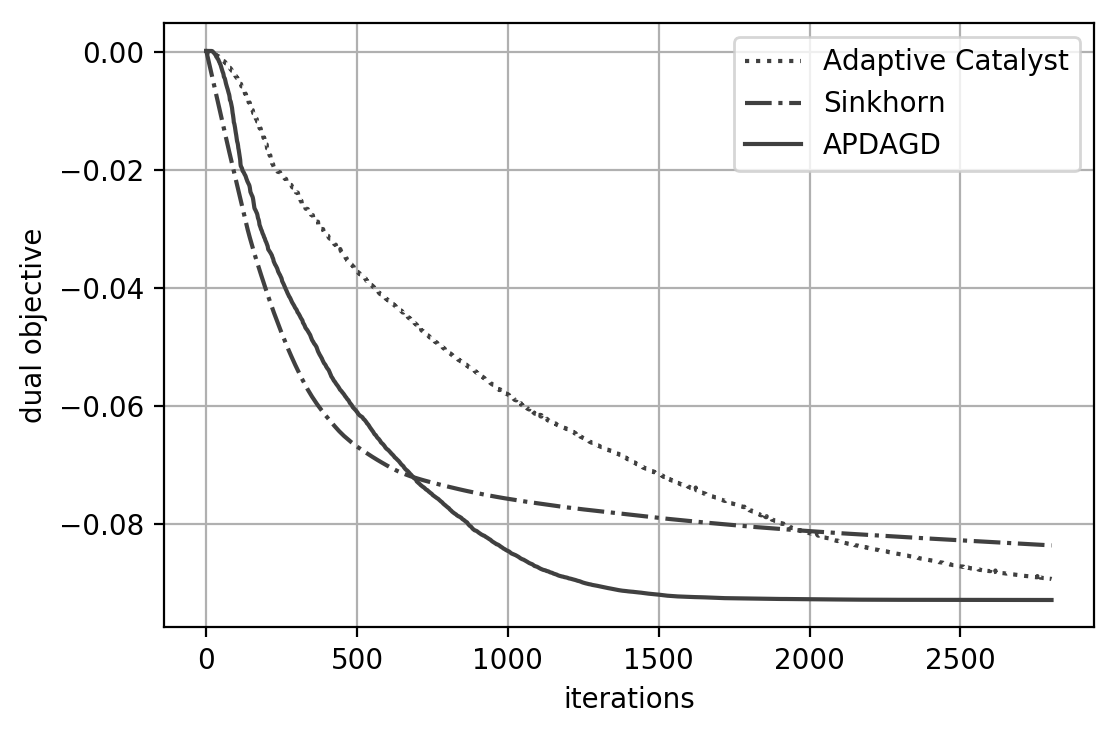}
\smallskip
\caption{Сходимость алгоритма Синхорна и каталист-оболчки, примененной к ЭОТ}
\label{cat}
\end{figure}

В связи с этим было произведено сравнение ускоренного с помощью адаптивной каталист-оболочки градиентного метода с адаптивным выбором шага \cite{nesterov2014universal, kamzolov2020universal} для ЭОТ с алгоритмом Синхорна и методом APDAGD. Результаты представлены на рисунке \ref{cat}.

\paragraph{Заключение}
{Представленно доказательство нового режима сходимости для ранее известного метода AGMsDR из \cite{nesterov2020primaldual} и для обобщения этого метода на задачи, допускающие альтернативную минимизацию из \cite{2019arXiv190603622G}. А именно режим линейной сходимости в случае задач для которых выполняется условие Поляка-Лоясиевича или сильно выпуклых задач, если константа $\mu$ в этих условиях неизвестна.

Произведено экспериментальное сравнение AGM и Alternating AGM и выялено расхождение с теоретическими оценками этих методов на примере задачи ЭОТ и задачи минимизации квадратичной функции. Представлено объяснение такого поведения, связанное с возможной адаптивностью к сильной выпуклости на блоках переменных метода Alternating AGM, которое требует более детального изучения.

Также представлено ускорение градиентного методоа с помощью адаптивной каталист-оболочки для задачи ЭОТ, и показана неоправданность применения данного алгоритма на практике для данной задачи.
}
\bibliography{references.bib}

\end{document}